\documentclass[preprint,12pt]{elsarticle} 
\usepackage[english]{babel}

\usepackage{enumerate}         
\usepackage{amsmath} 
\usepackage{amssymb}  
\usepackage{graphicx}      
\usepackage{epsfig}        
   
\usepackage{color}    
\usepackage{caption}

\usepackage{enumerate}   
\usepackage{amsmath}
\usepackage{enumitem}

\journal{XXX}  

\begin{document}

\makeatletter
\def\rank{\mathop{\operator@font rank}\nolimits}
\def\ln{\mathop{\operator@font ln}\nolimits}
\def\Re{\mathop{\operator@font Re}\nolimits}
\def\diag{\mathop{\operator@font Diag}\nolimits}
\newcommand\bolden[1]{{\boldmath\bfseries#1}}
\makeatother


\begin{abstract}
This paper deals with an analysis and design of robust, state-feedback control law uniform-asymptotically stabilizing at origin the system consisting of coupled $n$th--order ordinary differential equations  in the presence of a non-vanishing at $x=0$ or even unbounded on the time interval $[0,\infty)$ time-varying high-frequency oscillating perturbation $w(t,x).$ The obtained results generalize and extend some known and now classical results in the control theory for a wider class of perturbations. Moreover, as is shown in the paper, there is no room for further generalization for $w$ which is time-dependent only, $w=w(t).$
\end{abstract} 

\title{Asymptotic stabilization of a system of coupled $n$th--order differential equations with potentially unbounded high-frequency oscillating perturbations}
\author{R.~Vrabel}
\ead{robert.vrabel@stuba.sk}
\address{Slovak University of Technology in Bratislava, Institute of Applied Informatics, Automation and Mechatronics,  Bottova 25,  917 01 Trnava,   Slovakia}

\newtheorem{thm}{Theorem}
\newtheorem{lem}[thm]{Lemma}
\newtheorem{defi}[thm]{Definition}
\newdefinition{rmk}{Remark}
\newdefinition{ex}{Example}
\newproof{pf}{Proof}

\maketitle

\section{Introduction}

Suppose we are interested in assessing and ensuring the robustness of the nominal control system $\dot x=f(x,u)$ against modeling errors, system uncertainties, external disturbances, {\it etc.}, represented by the perturbation term $w(t,x)$ added to the right side of the nominal system, 
\begin{equation}\label{def_system}
\dot x=f(x,u)+w(t,x), \ t\geq0.
\end{equation}
The state vector $x,$ control input $u,$ vector field $f$ and perturbation term $w$ are the vectors of suitable dimensions, provisionally let $x\in\mathbb{R}^d,$ $u\in\mathbb{R}^m,$ $d\geq m\geq1.$ We always assume that $f$ and $w$ are at least continuous and that $f(0,0)=0.$ The perturbations $w$ are assumed to be potentially unknown  but belonging to the class of {\it diminishing} functions  which covers the high-frequency oscillating and among them also some unbounded perturbations (Definition~\ref{def:diminishing} and Remark~\ref{diminishing_w}). Further, let us assume that the solutions $x$ of (\ref{def_system}) for each admissible control $u$ are unique to the right, that is, $x(t; t_0, x_0)$ is uniquely determined by $(t_0, x_0)$ for $t\geq t_0\geq0$. 

For a motivation, let us consider $x=0$ being an uniform-asymptotically stable equilibrium point of the nominal system $\dot x=f(x,u)$ for a state-feedback control $u=g(x).$  What can we say about the stability of any kind for the perturbed system? This question represents one of the fundamental problems in the various areas of robust stabilization of the control systems, see e.~g. \cite{Bagherzadeh}, \cite{Li}, \cite{Ma}, \cite{Zhu},
and in principle, to answer this question, it makes usually a difference whether the origin remains an equilibrium for the perturbed system or not. If $w(t,0)=0$, then the origin is an equilibrium of (\ref{def_system}). In this case, then we can analyze the stability behavior of the origin as an equilibrium of the perturbed system. If $w(t,0)\neq0$, then the origin is no longer an equilibrium of (\ref{def_system}). In this case, we usually analyze the ultimate boundedness of the solutions of the perturbed system. As have been shown in \cite[Chapter~9]{Khalil} if for an appropriate choice of the control law $u(t,x)$ the point $x=0$ becomes an exponentially stable equilibrium point of the nominal system and the perturbation term $w$ satisfies
\begin{equation}\label{bound} 
|w(t,x)|\leq\gamma(t)|x|+\eta(t),\ \forall|x|<r,\ \forall t\geq0 
\end{equation}
where $\gamma,\eta: [0,\infty)\rightarrow[0,\infty)$ are continuous, $\int_0^{\infty}\gamma(\tau)d\tau<\infty$ and $\eta$ is bounded, then for $\eta\equiv0,$ the origin is an exponentially stable equilibrium point of perturbed system and the solutions of perturbed system are ultimately bounded in the opposite case (that is, if $\eta$ is not identically zero). These analyses are close to the notion input-to-state stability which has been introduced by E. Sontag in \cite{Sontag}. In contrast to the case of exponential stability, a nominal system with uniform-asymptotically stable (but not exponentially stable) origin is not robust to the smooth perturbations with arbitrarily small linear growth bounds of the form $|w(t,x)|\leq\gamma|x|,$ $|x|<r,$ $t\geq0$ and $\gamma>0,$  see \cite{Khalil} for more details. Definitions of the above concepts are given in the following section.   

\centerline{}

Summarizing these facts, the general framework for our considerations and analyses is that
\begin{itemize}
\item[1)] we will assume the stabilizability of the nominal system $\dot x=f(x,u)$ at $x=0$ by a continuously differentiable state-feedback control $u=g(x),$ $g(0)=0.$ This property is guaranteed by the non-singularity assumption of Jacobian matrix of the function $f$ with respect to the variable $u$ at the point $(0,0)$ (Theorem~\ref{thm1});
\item[2)] we will not assume that $w(t,x)$ satisfies the inequality constraint of the form (\ref{bound}) and therefore the classical results of  Khalil \cite[Lemma~9.4, p.~352]{Khalil} based on the Lyapunov's converse theorem, Coddington \& Levinson \cite[Theorem~3.1, p.~327]{Coddington_Levinson}, Hartman \cite[Chapter X]{Hartman} both based on the state-space model representation, and their various variants (e.g.~\cite{Brauer}, \cite{Ladde}, \cite[p.~183]{Brauer_Nohel}, \cite{Struble}) are not applicable here in general. Moreover, for $\eta(t)$ bounded, but non-vanishing at $t=\infty,$ we obtain the stronger result by considering the subclass of diminishing functions $w(t,x)$ (Remark~\ref{diminishing_w}, Part~\ref{item2}), namely, vanishing of $x(t;t_0,x_0)$ at $t=\infty$ versus boundedness of $x(t;t_0,x_0)$ only. We will also consider the perturbations $w(t)$ that are unbounded for $t\rightarrow \infty$ (Example~\ref{example}). Our approach come out from the impressive results and deep theory developed by Strauss \& Yorke in \cite{Strauss_Yorke}, whose results are obtained by a thorough and fine analysis of solution behavior.
\end{itemize}
It is a known fact that in the linear case, $f(x,u) = Ax+Bu,$ a necessary and sufficient condition for stabilization (by a linear feedback law $u=Gx$ and in the sense of the pole placement problem) is that rank of the controllability matrix $[B,AB,\dots,A^{d-1} B] = d,$ which in turn is equivalent to the complete controllability in the open-loop sense. The situation is quite different in the nonlinear case. The control system with dynamics   
\[
f(x_1, x_2, x_3, u_1, u_2) = (u_1, u_2\cos x_1, u_2\sin x_1) 
\]
is easily seen to completely controllable, however, the system cannot be stabilized to $0$ by a $C^1$ state-feedback because of the Brockett's necessary condition for feedback stabilization (\cite[p.~186]{Brockett}). To see that above example does not satisfy Brockett's necessary condition, note that the points $(b,0,a)$ with $a\neq0$ are not contained in the image of $f$ for $x_1$ near $0-$see also Remark~\ref{remark_2} at the end of the paper. This is another demonstration of subtlety of the concept "stabilizability" for nonlinear systems.
  
This paper has been written in order to provide theoretical background for extending the existing results regarding the {\it eventual} uniform-asymptotic stabilizability of control systems  at origin (Definition~\ref{def:EvUAS}) by a  continuous state-feedback control law $u=g(x)$  to a wider class of admissible perturbations $w,$ namely, when the perturbing term is diminishing (Theorem~\ref{thm1}). Moreover, as is shown in the mentioned theorem, there is no room for further generalization if $w(t,x)$ is time-dependent only, $w(t,x)\equiv w(t).$ There do not seem to be any results in the literature for control of the systems with this kind of perturbations, especially if these systems are affected by (potentially unbounded) high-frequency oscillation disturbance source. 
\section{Notations and definitions} 
Let 
\begin{itemize}
\item[] $\mathbb{R}^d$ denotes the finite-dimensional Euclidean $d-$space and let $|\cdot|$ denotes any $d-$dimensional norm and we will use $\Vert\cdot\Vert$ for the Euclidean norm. For later reference, recall that all norms on $\mathbb{R}^d$ are equivalent (\cite[p.~273]{HornJohnson}), so, $\theta_1|\cdot|\leq\Vert\cdot\Vert\leq\theta_2|\cdot|$ for some positive constants $\theta_1$ and $\theta_2$ depending on $|\cdot|$; 
\item[] $\Vert\cdot\Vert_{op}$ represents the operator norm induced by the norm $\vert\cdot\vert;$ 
\item[] $\diag(A)$ denotes the column vector of the main diagonal elements of the matrix $A;$ 
\item[] for $r>0$ and a fixed point $x_0\in\mathbb{R}^d,$  $B_r(x_0)\triangleq\left\{x\in\mathbb{R}^d:|x-x_0|<r\right\};$
\item[] the superscript $'\,T\,'$ is used to indicate transpose operator.
\end{itemize}
We now turn to the definitions of the stabilities, stated for the closed-loop system (\ref{def_system}) with a general (continuous) feedback $u=g(t,x),$ that we will use here, and that we have adopted and adapted from \cite{Strauss_Yorke}. 
\begin{defi}\label{def:EvUS} 
The origin is {\bf eventually uniformly stable} (EvUS) if
for every $\varepsilon>0$, there exists $\alpha = \alpha(\varepsilon)\geq0$ and $\delta =\delta(\varepsilon)>0$ such that
\[ 
|x(t; t_0, x_0)|<\varepsilon\ \ \mathrm{for\ all}\ \ |x_0|<\delta\ \ \mathrm{and}\ \ t\geq t_0\geq\alpha.
\]
It is {\bf uniformly stable} (US) if one can choose $\alpha(\varepsilon)=0.$
\end{defi}
\begin{defi}
The origin is {\bf eventually uniformly attracting} (EvUA) if there exist $\delta_0>0$ and $\alpha_0\geq 0$ and if  
for every $\varepsilon>0$ there exists $T=T(\varepsilon)\geq 0$ such that 
\[
|x(t; t_0, x_0)|<\varepsilon\ \ \mathrm{for\ all}\ \ |x_0|<\delta_0, \ \ t_0\geq\alpha_0,\ \ \mathrm{and}\ \  t\geq t_0+T. 
\]
It is {\bf uniformly attracting} (UA) if one can choose $\alpha_0=0.$
\end{defi}
\begin{defi}\label{def:EvUAS}  
The origin is {\bf eventually uniform-asymptotically stable} (EvUAS) if it is both EvUS and EvUA. It is {\bf uniform-asymptotically stable} (UAS) if it is both US and UA.
\end{defi}
As have been proved in \cite{Strauss_Yorke},  the concepts EvUAS and UAS are equivalent if and only if $x=0$ is a unique-to-the-right solution through $(t_0,0)$ of (\ref{def_system}) with $u=g(t,x)$ defined on $[t_0,\infty).$ So EvUAS  is a natural generalization of uniform-asymptotic stability in which it is not assumed that the zero function is a solution.  

These definitions (defined in  $\epsilon-\delta$ terms) are in fact equivalent to the following statements by using the special comparison functions known as class $\mathcal{K}$ and class $\mathcal{KL},$ \cite[p.~144 and also Lemma~4.5]{Khalil}:
\[
\mathrm{EvUS}\ \Leftrightarrow  (\exists \tilde\alpha\in\mathcal{K})\ |x(t; t_0, x_0)|\leq\tilde\alpha(|x_0|),\ \forall |x_0|<\delta, \ \forall t\geq t_0\geq\alpha
\]
and
\[
\mathrm{EvUAS}\ \Leftrightarrow (\exists \beta\in\mathcal{KL})\ |x(t; t_0, x_0)|\leq\beta(|x_0|,t-t_0),\ \forall |x_0|<\delta_0, \ \forall t\geq t_0\geq\alpha_0.
\]
The origin is globally EvUAS if and only if the last inequality is satisfied for any initial state $x_0$ ($\delta_0=\infty$).

A special case of UAS, so called exponential stability, arises when the class $\mathcal{KL}$ function $\beta$ takes the form $\beta(\tilde r,s)=\kappa \tilde re^{-\mu s},$ $\kappa, \mu>0.$
\begin{defi}
Let $w:\ [0,\infty)\times\mathbb{R}^d\rightarrow \mathbb{R}^d$ be continuous. Then $w$ is \bolden{vanishing at $ x=0$} if there exists $t^*\geq0$ such that for all $t\geq t^*$ is $w(t,0)=0;$ and $w$ is \bolden{vanishing at $t=\infty$} if there exists $r^*>0$ such that for all $|x|\leq r^*$ the function $w(t,x)\rightarrow 0$ for $t\rightarrow \infty.$
\end{defi}
\begin{defi}\cite{Strauss_Yorke}\label{def:diminishing} 
\ Let $h:\ [0,\infty)\rightarrow \mathbb{R}^d$ be continuous. Then $h$ is {\bf diminishing} if
\[
\sup\limits_{0\leq\lambda\leq1}\left\vert \int_t^{t+\lambda} h(\tau)d\tau\right\vert\rightarrow 0\ \mathrm{as}\ \ t\rightarrow \infty.
\]
\end{defi}
\begin{rmk}\label{diminishing_w}
\begin{enumerate}[label=(\textbf{P\arabic*})]

\item For example, if $h(t)\rightarrow 0$ as $t\rightarrow \infty$ then $h$ is diminishing. But, vanishing of $h(t)$ at $t=\infty$ is a sufficient condition only, not a necessary one. Indeed, let us consider 
\[
h(t)=\left(\cos (e^t),\sin (e^t),0,\dots,0\right).
\]
Then $h$ is diminishing; for any $\lambda\geq0,$ by integrating by parts with $a=e^{-\tau}$ and $\dot b=e^{\tau}\cos (e^{\tau}),$ we get
\[
\left\vert\int\limits_{t}^{t+\lambda}\cos (e^\tau) d\tau\right\vert\leq 2e^{-t}(1+e^{-\lambda})\leq 4e^{-t}.
\]
The same inequality holds for the second component of $h,$ and thus
\[
\left\Vert \int_t^{t+\lambda} h(\tau)d\tau\right\Vert\leq\sqrt{32}e^{-t}\rightarrow 0,
\] 
but $\Vert h(t)\Vert=1.$ The diminishing function may not be even bounded on $[0,\infty).$ For example,
\[
h(t)=\left(t\cos (t^4),t\sin (t^4),0,\dots,0\right)
\]
is diminishing as follows from the asymptotic properties of the Fresnel functions for large $t$, see, e.g., \cite[p.~149]{Bateman} or \cite{Weisstein}, and $\Vert h(t)\Vert=t\rightarrow \infty.$  In both cases, the functions $h$ represent high-frequency oscillations, bounded and unbounded, respectively.  These functions, depicted in Fig.~\ref{perturbation_h1}, will be used later in Example~\ref{example} to demonstrate the effectiveness of the proposed controller.
\begin{figure}[ht] 
\captionsetup{singlelinecheck=off}
   \centerline{
    \hbox{
     \psfig{file=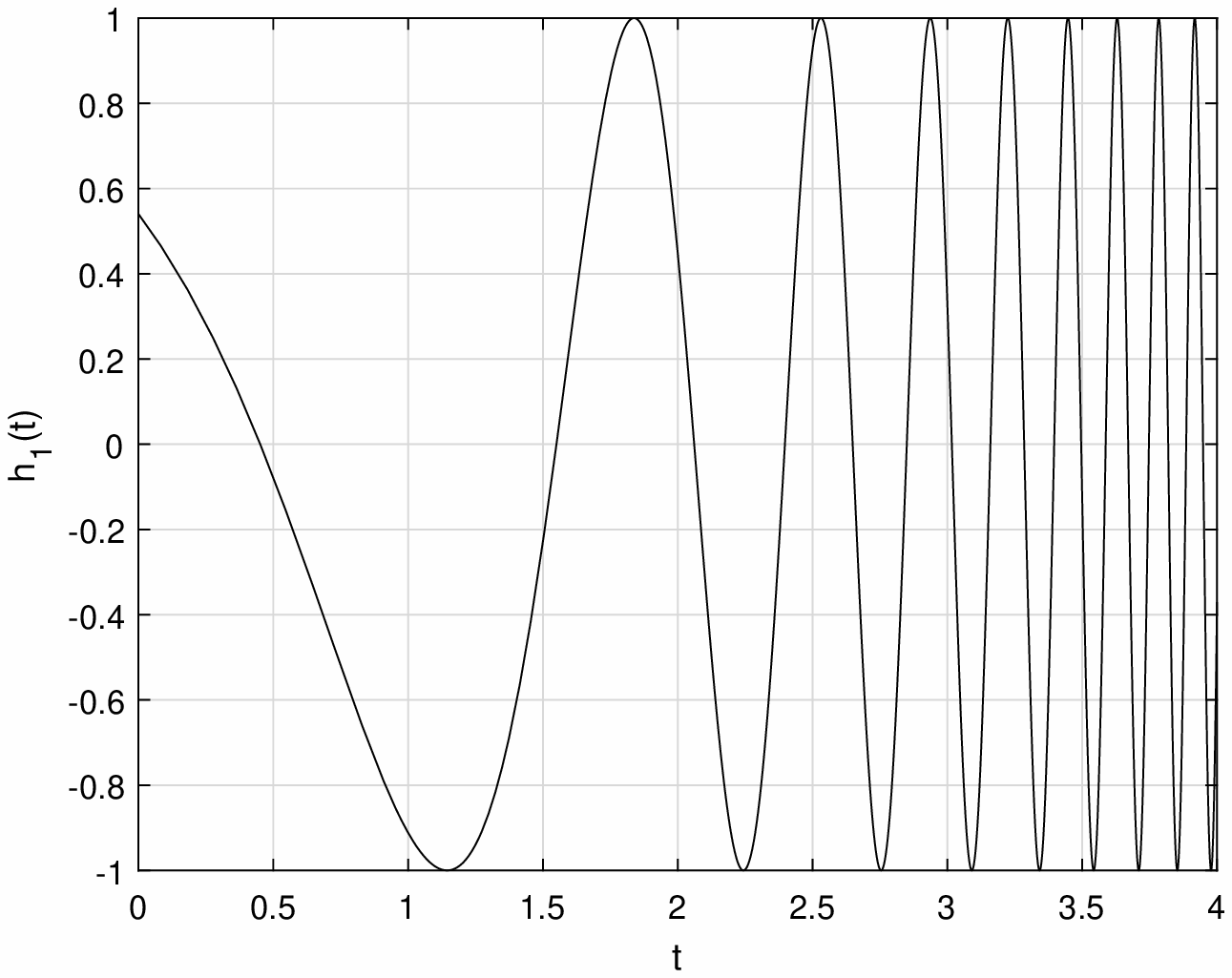,width=5.cm, clip=}
     \hspace{1.cm}
     \psfig{file=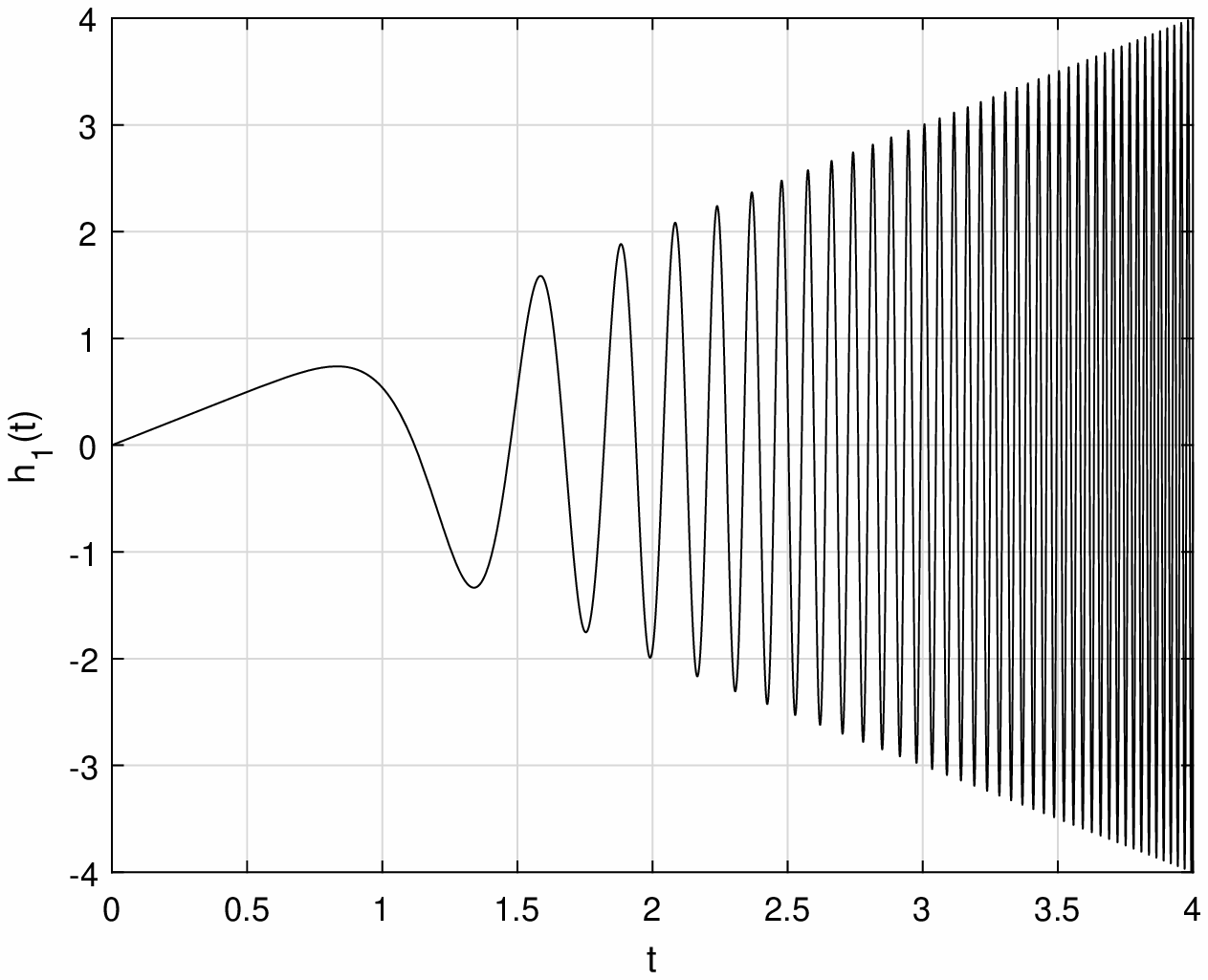,width=5cm,clip=}
    }
   }
\caption{ The functions $h_1(t)=\cos (e^t)$ (on the left) and $h_1(t)=t\cos (t^4)$ (on the right) on the interval $[0,4]$.}
\label{perturbation_h1}
\end{figure} 

\item \label{item2}  The concept of diminishing function can be naturally generalized for the functions depending also on $x,$ see \cite[Definition~2.19]{Strauss_Yorke} and the following discussion. We restrict ourselves to the diminishing functions of the form $w(t,x)=D(t)k(x),$ where each column of the $d\times d$ matrix $D$ is bounded on $[0,\infty)$ and diminishing in the sense of Definition~\ref{def:diminishing} and vector function $k: \mathbb{R}^d\rightarrow \mathbb{R}^d$ is continuous, that is, $k(x)\equiv1$ is also allowed. The boundedness of the columns of $D$ is required only if $k$ is a non constant function.
\end{enumerate}
\end{rmk} 
\section{Application to the system of n-order ODEs}  
In the framework given by the definitions above, our aim is to prove a new theorem on the eventually uniform-asymptotic stabilizability of the origin $0$ for the controlled system of $m$ $n\mathrm{th}-$order ordinary differential  equations ($m\geq1,$ $n>1$), which we may think as a special case of the system (\ref{def_system}),
\begin{equation}\label{def_system_n_order}
Y^{(n)}=F(Y,Y^{(1)},\dots,Y^{(n-1)}, U)+W(t,Y,Y^{(1)},\dots,Y^{(n-1)}),\ \ t\geq 0,
\end{equation}
given that $Y\in\mathbb{R}^{m},$ $U\in\mathbb{R}^{m},$ $F=(f_1,\dots,f_m)$ is $C^1$ function from $\mathbb{R}^{mn+m}$ to $\mathbb{R}^{m},$  the perturbation $W=(w_1,\dots,w_m)$ is continuous  from $[0,\infty)\times\mathbb{R}^{mn}$ to $\mathbb{R}^{m},$ and $Y^{(i)}$ denotes the $i-$th derivative with respect to the time $t,$ ($i=0,1,\dots,n-1$), and of course, we identify $Y^{(0)}$ with $Y$. For example, the Lagrange's equations in mechanics produce $m$ second--order differential equations for an $m-$degree of freedom dynamical system (\cite[Chapter~1]{Goldstein}, \cite[p.~158]{Murray_et_al}, \cite[p.~211]{Slotine_Li}, \cite[p.~435]{Vidyasagar}).

\begin{figure}[ht] 
\captionsetup{singlelinecheck=off}
   \centerline{
    \hbox{
     \psfig{file=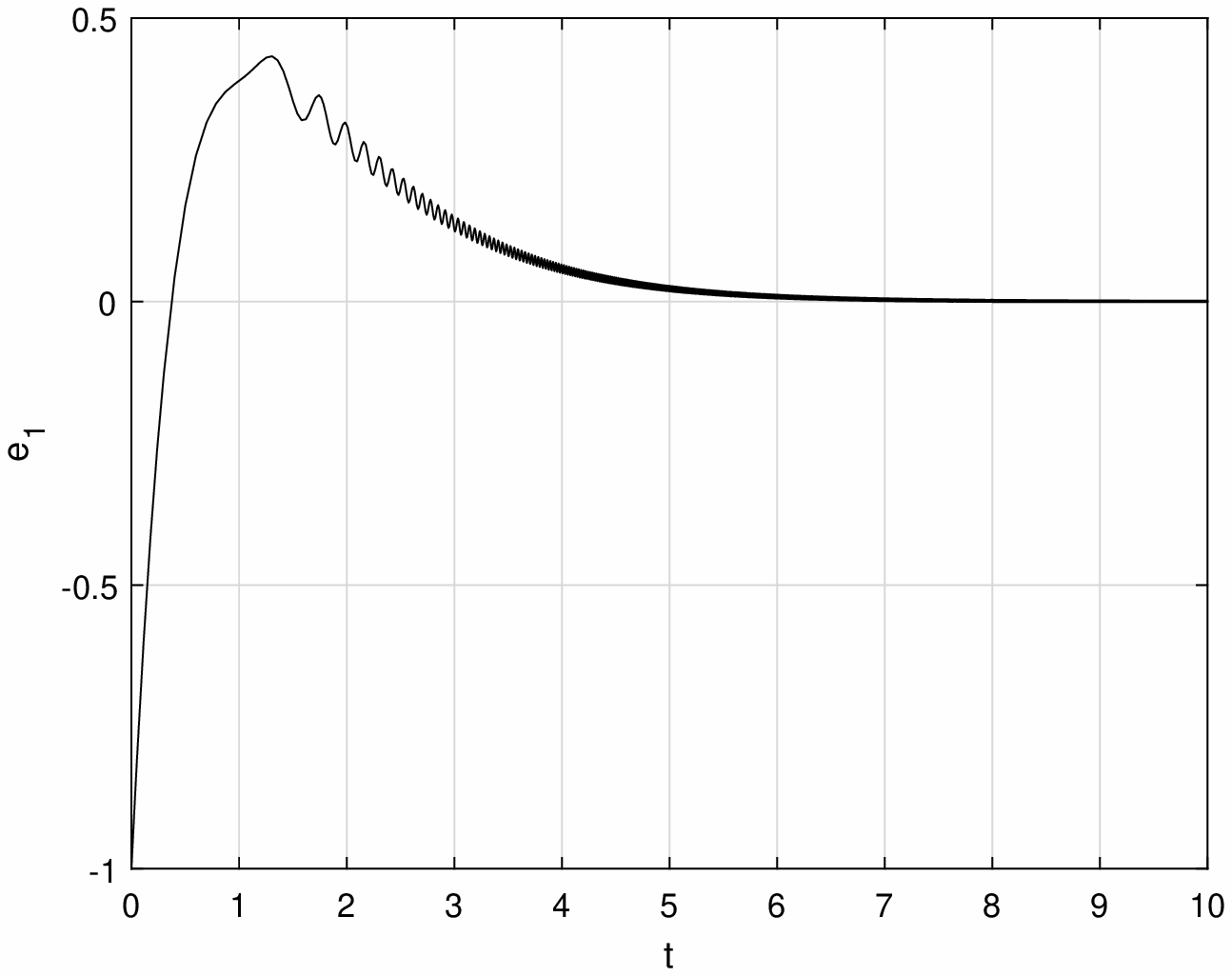,width=5.cm, clip=}
     \hspace{1.cm}
     \psfig{file=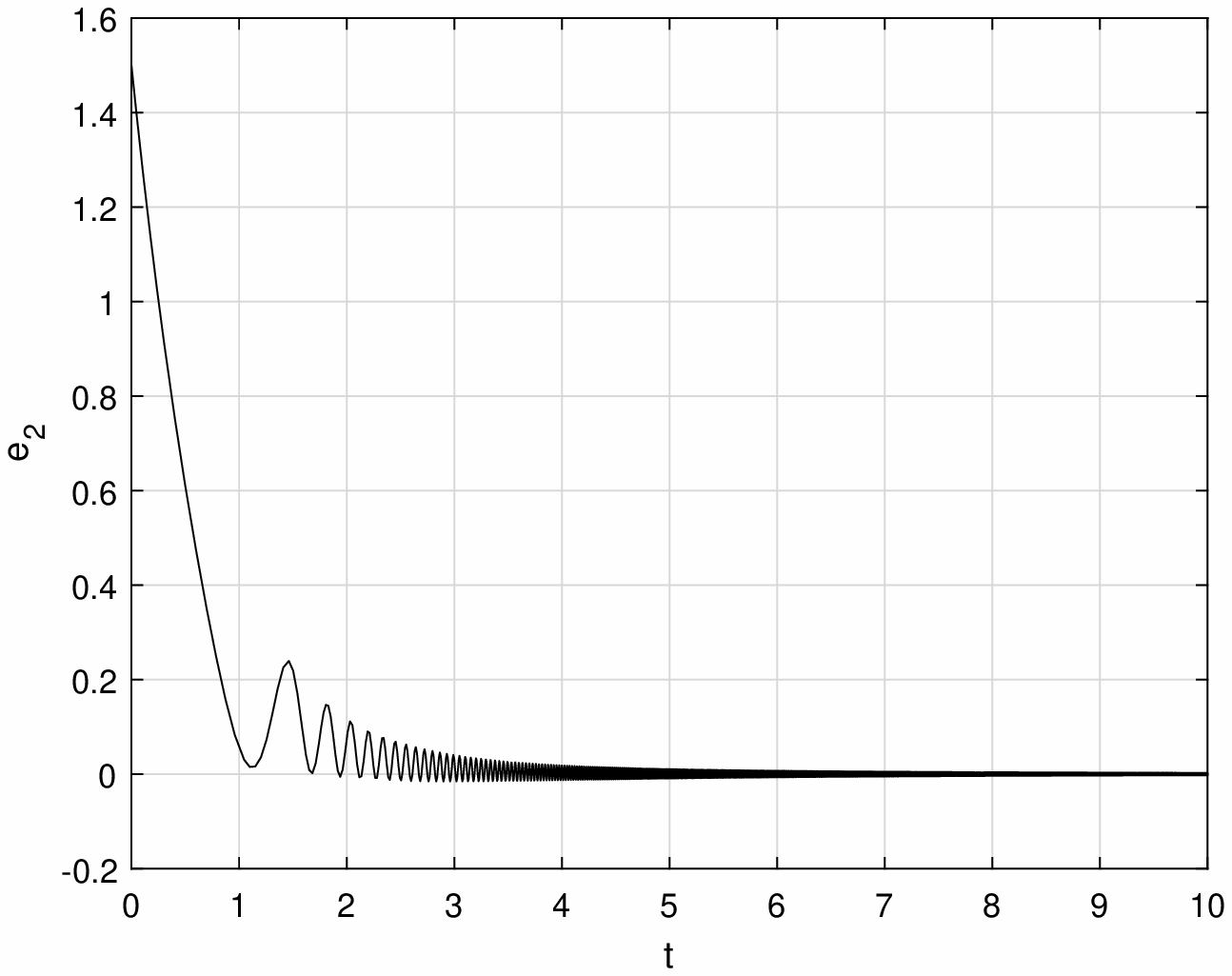,width=5cm,clip=}
    }
   }
   \vspace{0.5cm}
   \centerline{
    \hbox{
     \psfig{file=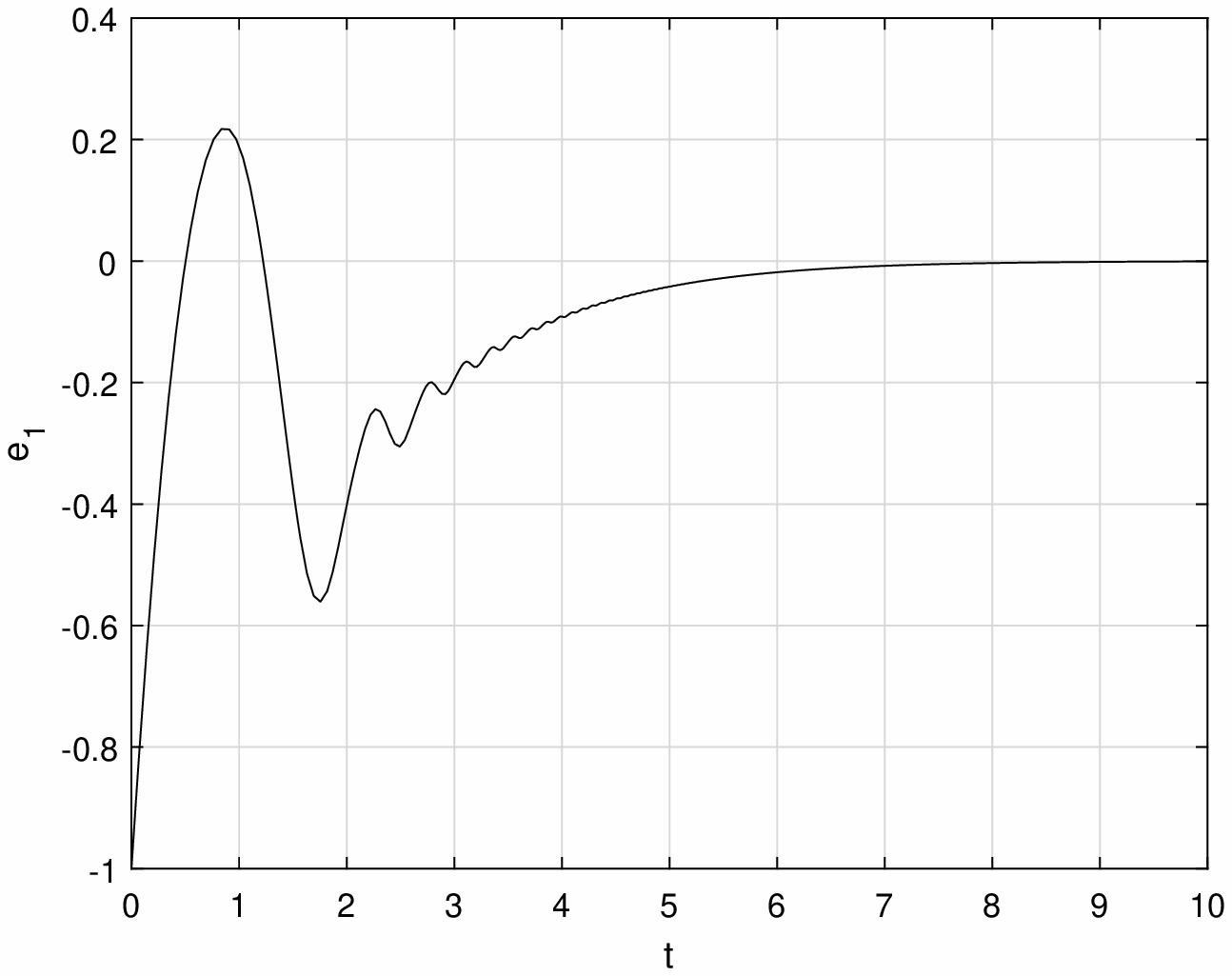,width=5.cm,clip=}
     \hspace{1.cm}
     \psfig{file=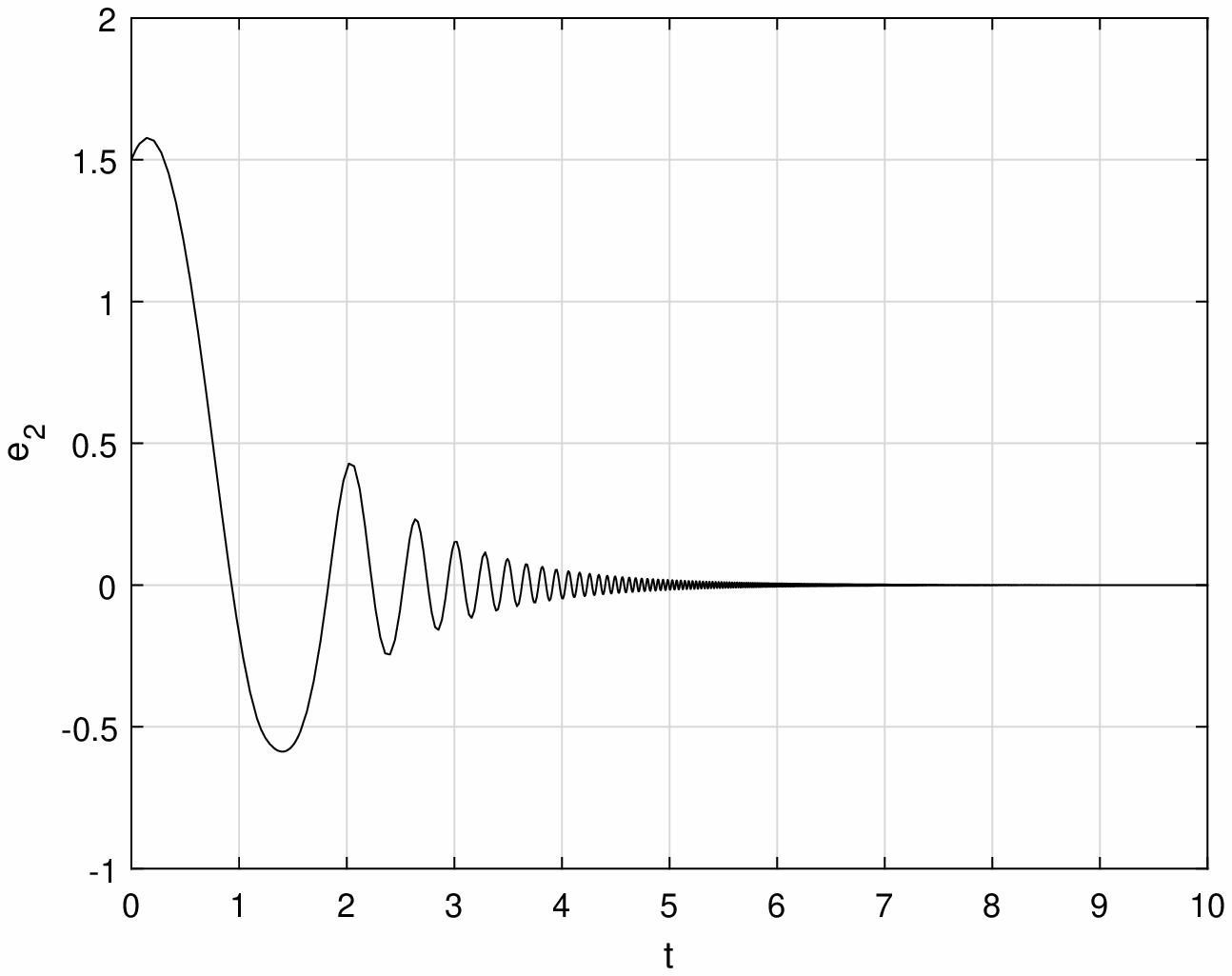,width=5.cm,clip=}
    }
   }
\caption{ Time evolution of the error dynamics $\dot E=A_HE+\tilde W(t,E)$ with the matrix $A_H=\left(\protect\begin{smallmatrix}-1&2\\0&-1.5\protect\end{smallmatrix}\right)$ and the perturbing term $\tilde W(t)=\left(0.5t\sin(t^4),-t\cos(t^4)\right)^T$ (the top row) and $\tilde W(t,e_1,e_2)=\left(-e_2\sin(e^t),2(e_1^{1/3}+e_2+1)\cos(e^t)\right)^T$ (the bottom row), respectively.}
\label{error_ei}
\end{figure}

Associating $Y$ with $X_1,$ and $Y^{(i)}$ with $X_{i+1},$ $i=1,\dots,n-1,$ we get the state variable matrix $X\triangleq[X_1,\dots,X_n]=[Y,Y^{(1)},\dots,Y^{(n-1)}]\in\mathbb{R}^{m\times n}$ and the system (\ref{def_system_n_order}) can be rewritten into the state-space representation 
\begin{equation}\label{state_space_model}
\dot X_i=X_{i+1}, \ i=1,\dots,n-1, \ \dot X_n=F(X,U)+W(t,X).
\end{equation}
Our main result is the following 
\begin{thm}\label{thm1}
Consider the control system (\ref{state_space_model}). Let $F:\ \mathbb{R}^{mn+m}\rightarrow \mathbb{R}^{m}$ is $C^1$ function, $F(0,0)=0$ and the corresponding Jacobian matrix with respect to the input variable vector $U$
\[
J_{F,U}(0,0)\triangleq\frac{\partial(f_1,\dots,f_m)}{\partial(u_1,\dots,u_m)}(0,0)
\]
is non-singular (and so bijective on $\mathbb{R}^{m}$). Let 
\begin{equation}\label{perturbation}
W(t,X)=D(t)K(X),
\end{equation}
where each column of an $m\times m$ matrix $D$ is bounded (for non constant $K(X)$) and diminishing, and the vector function $K=(k_1,\dots,k_m)^T:$ $\mathbb{R}^{mn}\rightarrow\mathbb{R}^m$ is continuous. 

Then there exists a $C^1$ state-feedback control law $U=G(X)$ defined in some open neighborhood of $0\in\mathbb{R}^{mn},$ $G(0)=0\in\mathbb{R}^m$ such that $0\in\mathbb{R}^{mn}$ is EvUAS for (\ref{state_space_model}) with $U=G(X).$ For $W$ independent of the state variable $X$, the EvUAS of $0$ for (\ref{state_space_model})  implies that $W(t)$ is diminishing. Moreover, if
\begin{itemize}
\item[a1)] for any $X\in\mathbb{R}^{mn}$ the functional $\Phi_X:\ \mathbb{R}^{m}\rightarrow \mathbb{R}$ given by the formula
\[
\Phi_X(U)=\frac12\vert F(X,U)\vert^2 
\]
is coercive, that is, $\lim_{\vert U\vert\rightarrow \infty}\Phi_X(U)=\infty,$
\item[a2)] the Jacobian matrix $J_{F,U}(X^*,U^*)$ is bijective for any $(X^*,U^*)\in\mathbb{R}^{mn+m},$
\end{itemize}
then the state-feedback control law $U=G(X)$ is defined globally, that is, for all $X\in\mathbb{R}^{mn}.$
\end{thm}
\begin{pf}
Let us define 
\begin{itemize}
\item the matrix $X_d(t)=\left[Y_d(t),Y_d^{(1)}(t),\dots,Y_d^{(n-1)}(t)\right]\in\mathbb{R}^{m\times n}$ that represents a desired trajectory of the control system (\ref{def_system_n_order}) -- in our case of asymptotic stabilization $X_d\equiv0,$ 
\item the matrix $\Delta=X-X_d(t),$ 
\end{itemize}
and
\begin{itemize}
\item let the tracking error is given as
\begin{equation}\label{tracking_error_eq}
E\triangleq\diag\left(\Delta
\begin{bmatrix}
\Gamma \\
1 
\end{bmatrix}
\right),\ E=(e_1,\dots,e_m)^T\in\mathbb{R}^m,\ 
\begin{bmatrix}
\Gamma \\
1 
\end{bmatrix}
\in\mathbb{R}^{n\times m},
\end{equation}
where each column of $\Gamma$ is such that the polynomials $\gamma_{1,j}+\gamma_{2,j}z+\dots+\gamma_{n-1,j}z^{n-2}+z^{n-1},$ $j=1,\dots,m,$ have the roots which are either negative or pairwise conjugate with negative real parts, therefore $\Delta(t)\rightarrow0$ if $E(t)\rightarrow0$ for $t\rightarrow\infty$; the $n-$tuple $(\gamma_{1,j},\dots,\gamma_{n-1,j},1)$ is the $j$th column of $\begin{bmatrix}
\Gamma \\
1 
\end{bmatrix}.$
\end{itemize} 
Differentiating (\ref{tracking_error_eq}) we obtain
\begin{equation*}
\dot E=\diag\left(\Delta
\begin{bmatrix}
0 \\
\Gamma
\end{bmatrix}\right)+F(X,U)-Y_d^{(n)}(t)+W(t,X)\ \ (\mathrm{here,}\ Y_d^{(n)}\equiv0)
\end{equation*}
and hence
\[
\dot E=A_HE+\underbrace{\diag\left(\Delta
\begin{bmatrix}
0 \\
\Gamma
\end{bmatrix}\right)+F(X,U)-A_HE}_{\triangleq\tilde F(X,U)}+W(t,X),
\]
for provisionally arbitrary $m\times m$ constant matrix $A_H.$ Because
\[
\diag\left(\Delta
\begin{bmatrix}
0 \\
\Gamma
\end{bmatrix}\right)-A_H \diag\left(\Delta
\begin{bmatrix}
\Gamma \\
1 
\end{bmatrix}\right), \ \Delta\equiv X
\]
is independent of $U,$  the Jacobian matrix $J_{\tilde F,U}(0,0)=J_{F,U}(0,0)$ and therefore is non-singular, also $\tilde F(0,0)=0.$ On the basis of the implicit function theorem, see, e.g. \cite[p.~136]{Shirali_Vasudeva}, there exists a neighborhood $P$ of $0\in\mathbb{R}^{mn},$  a neighborhood $Q$ of $0\in\mathbb{R}^{m}$ and a class $C^1$ function $G: P\rightarrow Q$ such that $G(0)=0$ and for all $(X_1,\dots,X_n)\in P$ is $\tilde F(X_1,\dots,X_n,G(X_1,\dots,X_n))=0.$ 

For a given fixed initial state $(t_0,X(t_0)),$ the mapping between the vector's $E(t)$ individual components and the rows of $X(t;t_0,X(t_0))$ is one-to-one, therefore $X$ can be expressed in terms of $E$ by the variation-of-parameter method, $X=X(E).$ The rest of the proof of the first part of the theorem (local stabilizability property) follows by applying the following lemma.
\begin{lem}
Consider the error dynamics 
\begin{equation}\label{error_system}
\dot E(t)=A_HE(t)+\tilde W(t,E),\ t\geq0, \ E\in\mathbb{R}^m,
\end{equation}
where all eigenvalues of the matrix $A_H\in\mathbb{R}^{m\times m}$ have negative real parts and $\tilde W(t,E)=W(t,X(E))=D(t)K(X(E)).$ Then $0$ is globally EvUAS for (\ref{error_system}).
\end{lem}
\begin{pf}
The statement of lemma follows from \cite[Corrolary~4.5 and 4.6, and Theorem~A(i)]{Strauss_Yorke} applied to (\ref{error_system}).
\end{pf}
We still need to ensure to be $X(t)\in P$ for all $t\geq t_0.$ Let $B_{r_{\max}}(0)$ is the maximal open ball in $P.$ From Definition~\ref{def:EvUS}, $|E(t;t_0,E(t_0))|<\varepsilon$ for $t\geq t_0\geq\alpha$ if 
\[
|E(t_0)|=\left\vert\diag\left(X(t_0)
\begin{bmatrix}
\Gamma \\
1 
\end{bmatrix}\right) \right\vert<\delta_E(\varepsilon),
\]
that is, for $|X(t_0)|<\delta^*_E$ for some $\delta^*_E=\delta^*_E(\Gamma,\varepsilon)>0,$ which may be calculated from the inequality
\[
|E(t_0)|\leq\frac1{\theta_1}\left\Vert\diag\left(X(t_0)
\begin{bmatrix}
\Gamma \\ 
1 
\end{bmatrix}\right) \right\Vert\leq\frac{\gamma^*\sqrt{m}}{\theta_1}\left\Vert X(t_0)\right\Vert\leq\frac{\gamma^*\theta_2\sqrt{m}}{\theta_1}\left\vert X(t_0)\right\vert,
\]
where $\gamma^*=\max\left\{\mathrm{the\ absolute\ value\ of}\ \gamma_{i,j},\, 1,\ i=1,\dots,n-1,\, j=1,\dots,m \right\},$ and so, $\delta^*_E=\frac{\theta_1}{\gamma^*\theta_2\sqrt{m}}\delta_E(\varepsilon).$ 

The sufficiently small $\varepsilon$ is chosen such that $|X(t)|<r_{\max}$ (for $t\geq t_0\geq\alpha$) by estimating solutions to the system
\[
\Bigg\{\diag\left(X
\begin{bmatrix}
\Gamma \\
1 
\end{bmatrix}
\right)=\Bigg\}\diag\left(\left[Y,Y^{(1)},\dots,Y^{(n-1)}\right]
\begin{bmatrix}
\Gamma \\
1 
\end{bmatrix}
\right)=E(t) ,\ |E(t)|<\varepsilon,
\]
which is for $E(t)\equiv0$ globally exponentially stable. Thus,  for the suitable constants $\kappa\geq1,$ $\mu_\Gamma>0$ and $t\geq t_0,$ we obtain that
\[
\vert X(t)\vert\leq\kappa\vert X(t_0)\vert e^{-\mu_\Gamma(t-t_0)} +\int\limits_{t_0}^t e^{-\mu_\Gamma(t-\tau)}\vert E(\tau)\vert d\tau
\]  
\[
\leq\kappa\vert X(t_0)\vert e^{-\mu_\Gamma(t-t_0)}+\frac{\varepsilon}{\mu_{\Gamma}}\left(1- e^{-\mu_\Gamma(t-t_0)}\right)\leq\kappa\vert X(t_0)\vert +\frac{\varepsilon}{\mu_{\Gamma}}.
\]
Hence, $\vert X(t)\vert<r_{\max}$ if $\vert X(t_0)\vert<\frac1{\kappa}\left(r_{\max}-\varepsilon/\mu_{\Gamma}\right)\triangleq\delta^*_X(\Gamma,\varepsilon)$ and $\varepsilon<r_{\max}\mu_{\Gamma}.$ So, $X(t)\in P$ for $t\geq t_0\geq\alpha$ if $\vert X(t_0)\vert<\min\{\delta^*_E,\delta^*_X\}.$

The second part of theorem, the global stabilization property, is the consequence of \cite[Theorem~1]{Galewski_Radulescu} and the fact that a linear mapping on $\mathbb{R}^{m}$ given by the matrix $A_H$ is a globally Lipschitz function in the sense of definition in \cite[Section~4]{Strauss_Yorke} with the Lipschitz constant $L=\Vert A_H\Vert_{op}$ on the whole $\mathbb{R}^{m}$ and $W(t,X)$ of the form (\ref{perturbation}) is globally diminishing with regard to the variable $X$ \cite[Definition~2.19]{Strauss_Yorke}, namely, the mentioned definitions hold for the open balls $B_r(0)$ with center $E=0\in\mathbb{R}^{m}$ and $X=0\in\mathbb{R}^{mn}$ and each radius $r>0,$ respectively. The proof of Theorem~\ref{thm1} is complete.
\end{pf}
\begin{ex}\label{example}
As an illustrative example, let us consider the error dynamics of the form
\[
\dot E=\begin{bmatrix}
-1 & 2 \\
 0 & -1.5 
\end{bmatrix}E+\tilde W(t,E), \ t\geq0
\]
with the diminishing perturbation term 
\[
\tilde W(t)=\left(0.5t\sin(t^4),-t\cos(t^4)\right)^T
\] 
and 
\[
\tilde W(t,e_1,e_2)=\left(-e_2\sin(e^t),2(e_1^{1/3}+e_2+1)\cos(e^t)\right)^T,
\] 
respectively. The time evolution of the error $E(t)=(e_1(t),e_2(t))^T$ with an initial error value $E(0)=(-1,\, 1.5)^T$ are depicted in Fig.~\ref{error_ei}. Recall, that these perturbations do not satisfy the inequality (\ref{bound}), $\tilde W$ is unbounded in the first case and does not meet the inequality $\left|\tilde W(t,E)\right|\leq\gamma|E|+\eta$ for any $\gamma,\eta>0$ in the neighborhood of $E=0$ due to the $e_1^{1/3}$ in the second one.
\end{ex}

The paper will end with three remarks.
\begin{rmk}\label{remark_2}
It is now a classical result that there exists a linear and continuous stabilizing control law for $\dot x=f(x,u)$  with $f(0,0)=0$ provided the unstable modes of the linearized system are controllable and there exists no stabilizing control law $u=Gx$ if the linearized system has an unstable mode which is uncontrollable. For the problem considered here, for $n>1$ the number of state variables ($mn$) is greater than the control inputs ($m$). But from the specific form of nominal part of the system (\ref{state_space_model}), $\dot x_i=\hat F_i(x_1,\dots,x_{mn},u_1,\dots,u_m),$ $i=1,\dots,mn,$ the Jacobian matrix $J_{\hat F,x}(0,0)$ is directly, after an appropriate rearranging of the rows, in the canonical controllability form (\cite[p.~283]{Antsaklis_Michel}). This fact together with a non-singularity of $J_{F,u}(0,0),$ allowing the transformation of input matrix to the required canonical form, ensures the controllability of the linear part of above system. Therefore, does not matter how are distributed the eigenvalues of $J_{\hat F,x}(0,0)$ in the complex plane, all eigenvalues are controllable. These findings point to an alternative approach to the local asymptotic stabilization of nominal system by a linear state-feedback control law $u=Gx,$ where $G$ is a suitable $m\times mn$ constant matrix ensuring the asymptotic stability of linear part of closed-loop system, $J_{\hat F,x}(0,0)+J_{\hat F,u}(0,0)G.$
\end{rmk}
\begin{rmk}
For the practical computations, especially for the large matrices, here may be useful the sufficient condition to be the Jacobian matrix $J_{F,U}(0,0)$ non-singular, given by the implication: {\it If the matrix $B=(b_{ij})\in\mathbb{R}^{m\times m}$ is strictly diagonally dominant, that is, $|b_{ii}|>\sum_{\substack{j=1 \\ j\neq i}}^m |b_{ij}|$ for all $i=1,\dots,m,$ then $B$ is  non-singular.} This result is known as the Levy-Desplanques theorem, \cite[p.~349]{HornJohnson}.  
\end{rmk}
\begin{rmk}
As indicated in the first lines of the proof of theorem, its basic idea can be used also for a state-trajectory tracking problem with the obvious modifications at some places in the proof under the assumption that $F(X_d(t),0)\equiv0$ for $t\geq t_0.$ The function 
\[
\tilde F(t,\Delta,U)=
\diag\left(\Delta
\begin{bmatrix}
0 \\
\Gamma
\end{bmatrix}\right)
+F\left(\Delta+X_d(t),U\right)-Y_d^{(n)}(t)
-A_H \diag\left(\Delta
\begin{bmatrix}
\Gamma \\
1 
\end{bmatrix}\right)
\]
and $U=G(t,\Delta)=G\left(t,X-X_d(t)\right),$ $t\geq t_0\geq0.$
\end{rmk}
\section*{Conclusions}
In this paper we solved the problem of stabilizability of the control systems consisting of the coupled $n\mathrm{th}-$order differential equations and affected by the high-frequency oscillating perturbations $w(t,x)$ belonging to the class of diminishing functions, not necessary bounded and vanishing at $t=\infty$ or/and at $x=0.$  Under easily verifiable assumptions given in Theorem~\ref{thm1}, we have shown that there exists an $C^1$  state-feedback control law preserving an uniform-asymptotic stability of the closed-loop equilibrium point $x=0$ of the nominal (unperturbed) system.

\end{document}